\documentclass{amsart}

\usepackage{amssymb,latexsym,amsmath,amsfonts,hhline,comment,mathrsfs,stmaryrd}
\usepackage{pdfsync}
\usepackage{amsthm}
\usepackage{graphicx}
\usepackage{wrapfig}
\usepackage{caption}
\usepackage{subcaption}
\usepackage{indentfirst}
\usepackage{enumerate}
\usepackage{tikz}

\usepackage{ifpdf}
\ifpdf \usepackage[colorlinks=true, citecolor=blue, linkcolor=blue, urlcolor=blue]{hyperref} \fi

\newtheorem{formula}{}[section]
\newtheorem{definition}[formula]{Definition}
\newtheorem{corollary}[formula]{Corollary}
\newtheorem{remark}[formula]{Remark}
\newtheorem{lemma}[formula]{Lemma}
\newtheorem{theorem}[formula]{Theorem}

\def\thrm{\begin{theorem}}
\def\thrml#1{\begin{theorem}\label{#1}}
\def\ethrm{\end{theorem}}
\def\rmrk{\begin{remark}}
\def\rmrkl#1{\begin{remark}\label{#1}}
\def\ermrk{\end{remark}}
\def\dfntn{\begin{definition}}
\def\dfntnl#1{\begin{definition}\label{#1}}
\def\edfntn{\end{definition}}
\def\nmrt{\begin{enumerate}}
\def\enmrt{\end{enumerate}}
\def\tm#1{\item[{\rm (#1)}]}
\def\qtn{\begin{equation}}
\def\qtnl#1{\begin{equation}\label{#1}}
\def\eqtn{\end{equation}}
\def\lmm{\begin{lemma}}
\def\lmml#1{\begin{lemma}\label{#1}}
\def\elmm{\end{lemma}}
\def\crllr{\begin{corollary}}
\def\crllrl#1{\begin{corollary}\label{#1}}
\def\ecrllr{\end{corollary}}
\def\css{\begin{cases}}
\def\ecss{\end{cases}}
					
\def\proof{\noindent{\bf Proof}.\ }
\def\eprf{\hfill$\square$}

\def\cX{{\mathcal X}}

\def\fK{{\frak K}}

\DeclareMathOperator{\WL}{WL}

\def\qaq{\quad\text{and}\quad}
\def\qoq{\quad\text{or}\quad}
\def\sbs#1#2{_{\scriptscriptstyle{#1,#2}}}
\def\sbo#1{_{\scriptscriptstyle{#1}}}
\def\nrm#1{\|#1\|}

\DeclareMathOperator{\dimwl}{dim_{\scriptscriptstyle WL}}

\makeatletter 
\def\@seccntformat#1{\csname the#1\endcsname. } 
\def\@biblabel#1{#1.}

\makeatother

\title[The WL-dimension of chordal bipartite graphs without $T_2$]{The Weisfeiler-Leman dimension of chordal bipartite graphs without bipartite claw}
\author{Ilia Ponomarenko}
\address{St.Petersburg Department of the Steklov Mathematical Institute, St.Petersburg, Russia}
\email{inp@pdmi.ras.ru}

\author{Grigory Ryabov}
\address{Sobolev Institute of Mathematics, Novosibirsk, Russia}
\address{Novosibirsk State University, Novosibirsk, Russia}
\email{gric2ryabov@gmail.com}

\thanks{The work is supported by the Russian Foundation for Basic Research (project 18-01-00752)}

\date{}

\begin{document}

\begin{abstract}
A graph  $X$ is said to be chordal bipartite if it is bipartite and contains no induced  cycle of length at least~$6$. It is proved that if $X$ does not contain bipartite claw as an induced subgraph, then the Weisfeiler-Leman dimension of $X$ is at most~$3$. The proof is based on the theory  of coherent configurations.


\end{abstract}

\maketitle

\section{Introduction}

The Weisfeiler-Leman dimension (WL-dimension, for short)  of a finite graph~$X$ can roughly be thought as the minimum number $\dimwl(X)$ of variables in a formula of a natural fragment of first-order logic, which is valid only for graphs isomorphic to~$X$; equivalently, the graph $X$ is identified by the $d$-dimensional Weisfeiler-Leman algorithm with $d=\dimwl(X)$ (for details, see~\cite{Grohe2017}). The WL-dimension of a class~$\fK$ of graphs is defined to be 
$$
\dimwl{\fK}=\min_{X\in\fK}\dimwl(X).
$$

Interest in the WL-dimension in recent years caused, in particular, by the fact that if $\dimwl(\fK)$ is bounded from above by a constant~$d$, then the graph isomorphism problem restricted to~$\fK$ is solved in polynomial time by the $d$-dimensional Weisfeiler-Leman algorithm.  The graphs $X$ with $\dimwl(X)=1$ have  completely been characterized in~\cite{KSS}  and independently in~\cite{AKRV2017}. However, such a characterization for the graphs of the WL-dimension greater than one seems to be hopeless~\cite{FKV}. Moreover, there exist infinitely many graphs with  arbitrarily large WL-dimension~\cite{CFI}. 

\begin{figure}[h]
	\includegraphics[scale=0.4]{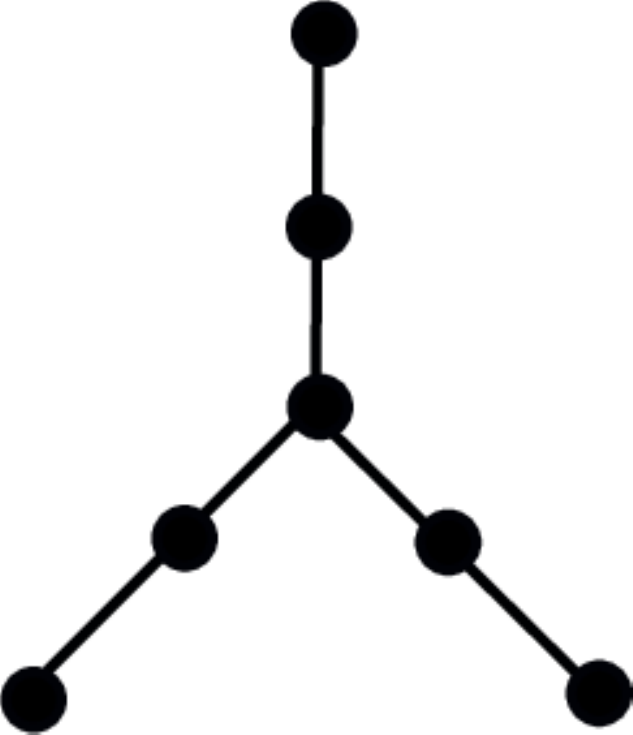}	
	\caption{The bipartite claw $T_2$.\label{fig1}}
\end{figure}

There are several results establishing an upper bound of $\dimwl{\fK}$ for specific classes~$\fK$, e.g., the planar graphs~\cite{KPS} or distance-hereditary graphs~\cite{GNP}.  In the present paper, we are interested in the WL-dimension of a special subclass of  chordal bipartite graphs (a bipartite graph is chordal if it contains no induced cycle of length at least~$6$).  The graph isomorphism problem for the class of all chordal bipartite graphs is polynomial-time equivalent to the graph isomorphism problem for general graphs~\cite{UTN}. Therefore, the WL-dimension of the chordal bipartite graphs is unlikely to be bounded from above by a constant. The subclass we mentioned consists of chordal bipartite graphs without bipartite claw, see Fig.~\ref{fig1}; a reason for this choice is that these graphs include several known classes (bipartite permutation graphs, difference graphs, etc.) with unknown WL-dimension. 

\thrml{main}
The $\WL$-dimension of the class of chordal bipartite $T_2$-free graphs  is equal to~$2$ or~$3$.
\ethrm

The proof of Theorem~\ref{main} is based on  theory of coherent configurations~\cite{CP}. A coherent configuration can be imagined as an arc-colored  complete graph with color classes satisfying some regularity conditions; these conditions are described via the so-called intersection numbers  (for exact definitions, see Sections~\ref{180520a} and~\ref{180520b}). 
According to~\cite{WeisL}, every graph $X$ is associated with uniquely determined coherent configuration~$\cX$ on the  vertex set of~$X$. It was proved in~\cite{FKV} that $\dimwl(X)\le 2$ if and only if $\cX$ is separable, i.e., is determined up to isomorphism by  the array of its intersection numbers. When one vertex of the graph $X$ is distinguished, the formula for $X$ from the definition of the WL-dimension should contain one more variable (this follows from~\cite{CFI}). Thus Theorem~\ref{main} is almost a direct consequence of the theorem below and the characterization of graphs with WL-dimension equal~$1$ (see above).

\thrml{main2} 
The coherent configuration of a connected chordal bipartite $T_2$-free graph with distinguished vertex is separable.
\ethrm

The proof of Theorems~\ref{main} and~\ref{main2} is given in Section~\ref{190520a}. The main tools for the proof are prepared in Section~\ref{08-520i}, where the coherent configurations of chordal bipartite graphs are studied. A relevant background on coherent configurations and graphs is given in Sections~\ref{180520a} and~\ref{180520b}, respectively.

\medskip

{\bf Notation.}

Throughout the paper, $\Omega$ is a finite set.

The set and number of classes of an equivalence relation $e$ on a set $\Omega$, are denoted by $\Omega/e$ and  $\nrm{e}=|\Omega/e|$, respectively.

For $r\subseteq \Omega\times\Omega$ and $\alpha\in\Omega$, we set  $\alpha r: =\{\beta\in \Omega:\ (\alpha,\beta)\in r\}$.

For $r\subseteq \Omega\times\Omega$ and $\Delta,\Gamma\subseteq \Omega$, we set $r\sbs{\Delta}{\Gamma}=r\cap(\Delta\times\Gamma)$ and put $r\sbo{\Delta}=r\sbs{\Delta}{\Delta}$.

The disjoint union of $m\ge 1$ copies of a complete bipartite graph with parts of cardinalities $a$ and $b$ is denoted by $mK_{a,b}$.

\section{Coherent configurations}\label{180520a}

In this section we provide a short background of the theory of coherent configurations. We use the notation and terminology from~\cite{CP}, where the most part of the material is contained. 

\subsection{Basic definitions} Let $\Omega$ be a finite set and $S$ a partition of $\Omega^2$; in particular, the elements of $S$ are treated as binary relations on~$\Omega$. A pair $\mathcal{X}=(\Omega,S)$ is called a \emph{coherent configuration} on $\Omega$ if the following conditions are satisfied:
\nmrt
\tm{C1}  the diagonal relation $\{(\alpha,\alpha):\ \alpha\in\Omega\}$ is the union of some relations of~$S$,
\tm{C2} for each $s\in S$, the relation $\{(\alpha,\beta): (\beta,\alpha)\in s\}$ belongs to~$S$,
\tm{C3}given $r,s,t\in S$ the number $c_{rs}^t=|\alpha r\cap \beta s^{*}|$ does not depend on $(\alpha,\beta)\in t$. 
\enmrt

Any relation belonging to $S$ is called a {\it basis relation} of~$\cX$. A set $\Delta \subseteq \Omega$ is called a \emph{fiber} of $\cX$ if the relation $\{(\delta,\delta):\ \delta\in\Delta\}$ is basis. The set of all fibers is denoted by $F=F(\cX)$. From the condition~(C1), it follows that
$$
\Omega=\bigcup_{\Delta\in F}\Delta
$$ 
and this union is disjoint. For each $s\in S$, there exist (uniquely determined) $\Delta,\Lambda\in F$ such that $s\subseteq\Delta\times\Gamma$. Moreover, if $\Delta$ or $\Gamma$ is a singleton, then $s=\Delta \times \Gamma$.\medskip

Let $\Delta$ be the union of some fibers of the coherent configuration~$\cX$. Denote by~$S_\Delta$ the set of all nonempty relations~ $s_\Delta$, $s\in S$. Then the pair $\cX_\Delta=(\Delta,S_\Delta)$ is a coherent configuration.

\subsection{Isomorphisms and separability} Let $\cX=(\Omega,S)$ and $\cX'=(\Omega',S')$ be two coherent configurations. A bijection $f:\Omega\to\Omega'$ is called a combinatorial isomorphism from $\cX$ to $\cX'$ if the relation 
$$
s^f=\{(\alpha^f,\beta^f):\ (\alpha,\beta)\in s\}
$$
belongs to $S'$ for every $s\in S$. The combinatorial isomorphism $f$ induces a natural bijection $\varphi:S\to S'$, $s\mapsto s^f$. One can see that $\varphi$  preserves the numbers from the condition~(C3), namely, the numbers $c_{rs}^t$ and $c_{r^{\varphi},s^{\varphi}}^{t^{\varphi}}$ are equal  for all $r,s,t\in S$. Every bijection $\varphi:S\to S'$ having this property is called an \emph{algebraic isomorphism} from $\cX$ to $\cX'$. A coherent configuration is called \emph{separable} if every algebraic isomorphism from it to another coherent configuration is induced by an isomorphism.

\subsection{Parabolics and twins}
An equivalence relation $e$ on the set $\Omega$ is called a {\it parabolic} of the coherent configuration~$\cX$ if $e$ is the union of some basis relations. Denote by $S_{\Omega/e}$ the set of relations 
$$
s\sbo{\Omega/e}=\{(\Delta,\Gamma)\in\Omega/e\times\Omega/e:\ s\sbs{\Delta}{\Gamma}\ne\varnothing\},\quad s\in S.
$$
Then the pair $\cX_{\Omega/e}=(\Omega/e,S_{\Omega/e})$ is a coherent configuration called the {\it quotient} of~$\cX$ modulo~$e$.\medskip

Following~\cite{GNP}, we say that  $\alpha,\beta\in \Omega$ are twins of $\cX$ or {\it  $\cX$-twins} if for each $\gamma\in\Omega$ other than $\alpha$ and $\beta$, and each $s\in S$, we have
$$
(\alpha,\gamma)\in s\quad\Leftrightarrow\quad (\beta,\gamma)\in s.
$$
It immediately follows that $\alpha$ and $\beta$ belong to the same fiber of~$\cX$. Moreover, the relation $t=t(\cX)$  ``to be $\cX$-twins'' is a parabolic of~$\cX$ \cite[Lemma~3.1]{GNP}; it is called a {\it twin parabolic} of the coherent configuration~$\cX$.

\lmml{160520a}\cite[Proposition 4.10]{GNP}
A coherent configuration $\cX$ is separable if the quotient of $\cX$ modulo $t(\cX)$ is separable.
\elmm

\subsection{Coherent closure}
There is a natural partial order\, $\le$\, on the set of all coherent configurations on the same set~$\Omega$. Namely, given two such coherent configurations $\cX$ and $\cX'$, we set $\cX\le\cX'$ if and only if each basis relation of~$\cX$ is the union union of some basis relations of~$\cX'$. The {\it coherent closure} $\WL(T)$ of a set $T$ of binary relations on $\Omega$, is defined to be the smallest coherent configuration on $\Omega$, for which each relation of~$T$ is the union of some basis relations. 

\section{Graphs}\label{180520b}

\subsection{Basic notation}
By a {\it graph} we mean a (finite) simple undirected graph, i.e., a pair $X=(\Omega,D)$ of a set $\Omega$ of vertices and an irreflexive symmetric relation $D\subseteq \Omega\times \Omega$, which represents the edge set of $X$. The elements of $D$ are called {\it arcs}, and $D$ is the {\it arc set} of~$X$. Two vertices $\alpha,\beta\in \Omega$ are said to be adjacent in $X$ whenever $(\alpha,\beta)\in D$; we also say that $\beta$ is the neighbor of $\alpha$ in~$X$. The subgraph of $X$ induced by $\Delta\subseteq \Omega$ is denoted by~$X_\Delta$.\medskip 

The graph $X$ is said to be {\it empty} if $D=\varnothing$.  A bipartite graph $X$ with parts~$\Delta$ and $\Gamma$ is said to be {\it biregular} if the number of $X$-neighbors of a vertex $\alpha\in\Omega$ depends only on whether $\alpha\in\Delta$ or $\alpha\in\Gamma$.\medskip 

The {\it distance} $d(\alpha,\beta)$ between the vertices $\alpha,\beta$ of $X$ is defined as usual to be the length of a shortest path in~$X$ from one of $\alpha,\beta$ to the other. The minimal distance from $\alpha$ to a vertex belonging to a set $\Delta\subseteq\Omega$ is denoted by $d(\alpha,\Delta)$. 

\subsection{Coherent configuration of a graph}
The coherent configuration $\WL(X)$ of  the graph $X$ is defined to be the coherent closure $\WL(\{D\})$.  In the lemma below we establish two properties of $\WL(X)$ to be used in the sequel.

\lmml{regular}
Let $X$ be a graph with arc set $D$, $\cX\geq \WL(X)$ a coherent configuration, and $\Delta,\Gamma\in F(\cX)$. Then 
\nmrt
\tm{1} the bipartite graph with parts $\Delta$ and $\Gamma$ and arc set $D\sbs{\Delta}{\Gamma}\cup D\sbs{\Gamma}{\Delta}$ is biregular,
\tm{2} if $\Gamma=\{\alpha\}$ for some $\alpha\in\Omega$, then $d(\alpha,\Delta)=d(\alpha,\delta)$ for all $\delta\in\Delta$.
\enmrt
\elmm 
\proof Statement~(1) follows from \cite[formula~(2.1.4)]{CP}. To prove statement~(2), denote by $s$ the set of all  pairs of vertices of $X$ at distance $d\ge 0$. According to \cite[Theorem~2.6.7]{CP}, $s$ is a union of basis relations of $\WL(X)$ and hence of~$\cX$. Now let $d=d(\alpha,\delta)$ for some $\delta\in\Delta$. Then $s\sbs{\Delta}{\Gamma}\ne\varnothing$ and hence is a basis relation of~$\cX$. Therefore, $s\sbs{\Delta}{\Gamma}=\Delta\times\{\alpha\}$ and we are done by the definition of~$s$.\eprf\medskip

Sometimes, it is convenient to consider a graph $X$ in which a certain vertex~$\alpha$ is fixed. In this case, we use notation~$X_\alpha$ and say that $X_\alpha$ is a graph with distinguished vertex~$\alpha$. The coherent configuration of~$X_{\alpha}$ is defined to be the coherent closure 
$$
\WL(X_\alpha)=\WL(\{D,1_\alpha\}),
$$ 
where $1_\alpha=\{(\alpha,\alpha)\}$.  One can see that $\WL(X_\alpha)\ge\WL(X)$ and $\{\alpha\}\in F(\cX)$ for any $\cX\ge\WL(X)$.

\lmml{bifiber}
Let $X$ be a connected bipartite graph with distinguished vertex and $\cX\ge WL(X)$. Then each $\Delta\in F(\cX)$ is contained in one of the two parts of~$X$. In particular, the graph $X_\Delta$ is empty.
\elmm

\proof Denote by $\alpha$ the distinguished vertex of $X$. The connectivity assumption implies that $d_X(\alpha,\beta)\ne\infty$ for all vertices $\beta$. Since the graph $X$ is also bipartite, the vertices $\alpha$ and $\beta$ belong to the same part of~$X$ if and only if the number $d_X(\alpha,\beta)$ is even. Thus the required statement follows from Lemma~\ref{regular}(2).\eprf

\subsection{Twins in graphs}
The vertices $\alpha$ and $\beta$ of the graph $X$ are called {\it twins} or {\it $X$-twins}  if any other vertex is the neighbor of both $\alpha $ and $\beta$ or none of them. One can see that the relation $e=e(X)$ ``to be $X$-twins'' is an equivalence relation on the vertex set of~$X$. 

 \lmml{twin1} \cite[Lemma~3.4(2)]{GNP}
Let $X$ be a graph with distinguished vertex. Then  every two $X$-twins belonging to the same fiber of $\WL(X)$ are also $\WL(X)$-twins.
\elmm

In the condition of Lemma~\ref{twin1}, let $\Delta,\Gamma\in F(\cX)$. We define  an equivalence relation $e\sbs{\Delta}{\Gamma}$ consisting of all pairs $(\alpha,\beta)\in\Delta\times\Delta$ such that $\alpha$ and $\beta$ are twins of the bipartite graph $X\sbs{\Delta}{\Gamma}$ with parts $\Delta$ and $\Gamma$ and the arc set $D\sbs{\Delta}{\Gamma}\cup D\sbs{\Gamma}{\Delta}$. It should be noted that $e\sbs{\Delta}{\Gamma}$ depends on $X$ but we avoid write $e(X\sbs{\Delta}{\Gamma})$ if the graph~$X$ is clear from the context.

\lmml{060520a}
In the above notation, the following statements hold,
\nmrt
\tm{1} $\nrm{e\sbs{\Delta}{\Gamma}}=1$ if and only if $D\sbs{\Delta}{\Gamma}=\Delta\times\Gamma$ or $D\sbs{\Delta}{\Gamma}=\varnothing$,
\tm{2} if $\Lambda\in F(\cX)$ and $e\sbs{\Delta}{\Gamma}\subseteq e\sbs{\Delta}{\Lambda}$, then $\nrm{e\sbs{\Delta}{\Gamma}}$ is divided by $\nrm{e\sbs{\Delta}{\Lambda}}$.
\enmrt
\elmm
\proof Statement~(1) is obvious. To prove statement~(2), we recall that $\WL(X\sbs{\Delta}{\Gamma})$ is the smallest coherent configuration containing the arc set of the graph $X\sbs{\Delta}{\Gamma}$. Since this set is the union of basis relations of the coherent configuration $\WL(X)\sbo{\Delta\cup\Gamma}$, we conclude that
$$
\WL(X\sbs{\Delta}{\Gamma})\le \WL(X)\sbo{\Delta\cup\Gamma}=\cX\sbo{\Delta\cup\Gamma}.
$$
However, from Lemma~\ref{twin1}, it follows that $e\sbs{\Delta}{\Gamma}=t(\WL(X\sbs{\Delta}{\Gamma}))\sbo{\Delta}$. By the above inclusion, this implies that $e\sbs{\Delta}{\Gamma}$ is a parabolic of $\cX_\Delta$. Similarly, one can verify that
$e\sbs{\Delta}{\Lambda}$ is also a parabolic of $\cX_\Delta$. Thus, the required statement follows from~\cite[Corollary 2.1.23, Exercise~3.7.9]{CP}.\eprf

\subsection{The Weisfeiler-Leman dimension}

The exact definition of the WL-dimen\-sion $\dimwl(X)$ of a graph $X$ requires a discussion about the multidimensional Weisfeiler-Leman algorithm, which is beyond the scope of the present paper; we refer the interested reader to the monograph~\cite{Grohe2017}. A relevant information on this invariant is collected in the lemma below.


\lmml{290420b}
Given a graph $X$, the following statements hold:
\nmrt
\tm{1}  $\dimwl(X)$ equals the maximum WL-dimension of a component of~$X$,
\tm{2}	$\dimwl(X)\le \dimwl(X_\alpha)+1$ for each vertex $\alpha$ of $X$, 
\tm{3}	$\dimwl(X)\le 2$ if and only if the coherent configuration $\WL(X)$ is separable. 
\enmrt
\elmm 
\proof  Statement~(1) is proved in the first paragraph of \cite[Section~4]{KPS}. Statement~(2) is an easy consequence of~\cite[Theorem 5.2]{CFI}. Statement~(3) follows from \cite[Theorem~2.1]{FKV}.\eprf

\section{Chordal bipartite graphs and coherent configurations}\label{08-520i}

A graph $X$ is said to be {\it chordal bipartite} if it is bipartite and contains no induced cycle of length at least~$6$, see~\cite[Section~12.4]{Gol}. One can see that every induced subgraph of $X$ is also chordal bipartite. According to \cite[Theorems~12.5, 12.8]{Gol}, each connected component of a nonempty chordal bipartite graph $X$ has at least two two {\it bisimplicial} vertices; by definition the vertices  $\alpha$ and~$\beta$ are bisimplicial if they are adjacent and $X_{\alpha D,\beta D}=X_{\alpha D\,\cup\, \beta D}$ is a complete bipartite graph (with parts $\alpha D$ and $\beta D$). \medskip

In this section, we are interested in how the fibers of the coherent configuration $\WL(X)$ dissect the arcs of~$X$. A key point here is the following easy observation. 

\lmml{bireg} 
Every nonempty biregular chordal bipartite graph is isomorphic to $mK_{a,b}$ for some integers $m,a,b\ge 1$.
\elmm
\proof Let $X$ be a nonempty chordal bipartite graph, and let $Y$ be a component of $X$ with parts $\Delta$ and $\Gamma$. 
Then there are bisimplicial vertices $\alpha\in\Delta$ and $\beta\in\Gamma$. Assume that $X$ is biregular.  Then
$$
|\alpha D|=|\alpha' D|\qaq |\beta D|=|\beta' D|
$$
for all $\alpha'\in\Delta$ and $\beta'\in\Gamma$, where $D$ is the arc set of~$X$. Since the graph~$X_{\alpha D\,\cup\, \beta D}$  is complete bipartite with parts $\alpha D$ and $\beta D$, this implies that
$$
Y=X_{\alpha D\,\cup\, \beta D}=K_{a,b},
$$
where $a=|\alpha D|$ and $b=|\beta D|$.  Since $X$ is biregular, these numbers do not depend on~$Y$. Thus, each connected component of $X$ is isomorphic to $K_{a,b}$, and the required  statement is true for $m$ being the number of connected components of~$X$.\eprf

\thrml{080520p}
Let $X$ be a chordal bipartite graph with distinguished vertex, and let $\cX=\WL(X)$ and $F=F(\cX)$. Then for each $\Delta\in F$, 
\nmrt
\tm{1} the graph $X_\Delta$ is empty,
\tm{2} if $\Gamma\in F$ and the graph $X_{\Delta\cup\Gamma}$ is not empty, then $X_{\Delta\cup\Gamma}=mK_{a,b}$, where $m=\nrm{e\sbs{\Delta}{\Gamma}}$ and $a,b\ge 1$,
\tm{3} $\nrm{e\sbs{\Delta}{\Gamma}}=\nrm{e\sbs{\Gamma}{\Delta}}$ for all $\Gamma\in F$,
\tm{4} the set $\{e\sbs{\Delta}{\Gamma}:\ \Gamma\in F\}$ is linear ordered with respect to inclusion.
\enmrt
\ethrm
\proof Statement~(1) follows from  Lemma~\ref{bifiber}. Let $\Gamma\in F$. Since $X_{\Delta\cup\Gamma}$  is an induced subgraph of~$X$, it is chordal bipartite. Moreover, $X_{\Delta\cup\Gamma}=X_{\Delta,\Gamma}$ is biregular by Lemma~\ref{regular}(1). By Lemma~\ref{bireg}, this implies that if  $X_{\Delta\cup\Gamma}$  is not empty, then it is isomorphic to $mK_{a,b}$ for some integers $m,a,b\ge 1$.  It follows that in this case, each class of $ e\sbs{\Delta}{\Gamma}$ (respectively, $e\sbs{\Gamma}{\Delta}$) is the intersection of $\Delta$ (respectively,~$\Gamma$) with  vertex set of a component of  $X_{\Delta\cup\Gamma}$. This proves statement~(2) and also statement~(3) except for the case when if $X_{\Delta\cup\Gamma}$  is empty. However, in the latter case, statement~(3) holds trivially.\medskip

To prove statement~(4), it suffices to verify that if $\Gamma,\Lambda\in F$, then
$$
e\sbs{\Delta}{\Gamma}\subseteq e\sbs{\Delta}{\Lambda}\qoq
e\sbs{\Delta}{\Lambda}\subseteq e\sbs{\Delta}{\Gamma}.
$$
Without loss of generality, we may assume that each of the graphs $X_{\Delta\cup\Gamma}$ and  $X_{\Delta\cup\Lambda}$ is nonempty. Then $X_{\Delta\cup\Gamma\cup\Lambda}$ being a chordal bipartite graph with parts $\Delta$ and $\Gamma\cup\Lambda$, contains bisimplicial vertices $\alpha\in\Delta$ and $\beta\in \Gamma\cup\Lambda$; for the definiteness, let $\beta\in\Gamma$. By the assumption, $\alpha$ has a neighbor $\lambda\in\Lambda$. In view of statement~(2),  this implies that
$$
\beta D\sbs{\Delta}{\Gamma}\in \Delta /e\sbs{\Delta}{\Gamma}\qaq \lambda D\sbs{\Delta}{\Lambda} \in \Delta/ e\sbs{\Delta}{\Lambda},
$$
where $D$ is the arc set of the graph~$X$. On the other hand, since the vertices $\alpha$ and~$\beta$ are bisimplicial, every neighbor of $\beta$ is adjacent to every neighbor of~$\lambda$. Thus,
$$
\beta D\sbs{\Delta}{\Gamma}\subseteq\lambda D\sbs{\Delta}{\Lambda}.
$$
This means that at least one class of $e\sbs{\Delta}{\Gamma}$ is contained in some class of $ e\sbs{\Delta}{\Lambda}$. Since $e\sbs{\Delta}{\Gamma}$ and $ e\sbs{\Delta}{\Lambda}$ are the parabolics of the coherent configuration $\cX_\Delta$ (having a unique fiber, namely~$\Delta$), this is possible only if $e\sbs{\Delta}{\Gamma}\subseteq e\sbs{\Delta}{\Lambda}$.\eprf

\section{Coherent configurations of  chordal bipartite $T_2$-free graphs}
 
The proof of Theorem~\ref{main2} is based on the following statement which refines Theorem~\ref{080520p}  for the chordal bipartite $T_2$-free graphs.

\thrml{300420a}
Let $X$ be a connected   chordal bipartite $T_2$-free graph with distinguished vertex, and let $\cX\ge\WL(X)$ be a coherent configuration. Then for all $\Delta,\Gamma\in F(\cX)$, 
\qtnl{010520l}
\nrm{e\sbs{\Delta}{\Gamma}}\le 2.
\eqtn
\ethrm
\proof In what follows, we set $F=F(\cX)$.  We need an auxiliary lemma.

\lmml{270420a8}
Let $\Delta,\Gamma,\Lambda\in F$ be pairwise distinct. Assume that $X_{\Delta\cup\Gamma}$ is not empty and  $e\sbs{\Delta}{\Gamma}\subsetneq e\sbs{\Delta}{\Lambda}$. Then
$$
 \nrm{ e\sbs{\Delta}{\Gamma}}=2\nrm{e\sbs{\Delta}{\Lambda}}.
$$
\elmm
\proof Suppose on the contrary that $\nrm{ e\sbs{\Delta}{\Gamma}}\,\ne 2\nrm{e\sbs{\Delta}{\Lambda}}$. By the lemma hypothesis, we have $\nrm{ e\sbs{\Delta}{\Gamma}}\ne \nrm{e\sbs{\Delta}{\Lambda}}$. Furthermore, $\nrm{ e\sbs{\Delta}{\Gamma}}$ is divided by~$\nrm{e\sbs{\Delta}{\Lambda}}$ by Lemma~\ref{060520a}(2). Thus, 
\qtnl{030520i}
\nrm{ e\sbs{\Delta}{\Gamma}}\ge 3\nrm{e\sbs{\Delta}{\Lambda}}.
\eqtn

\begin{figure}[h]
\includegraphics[scale=0.2]{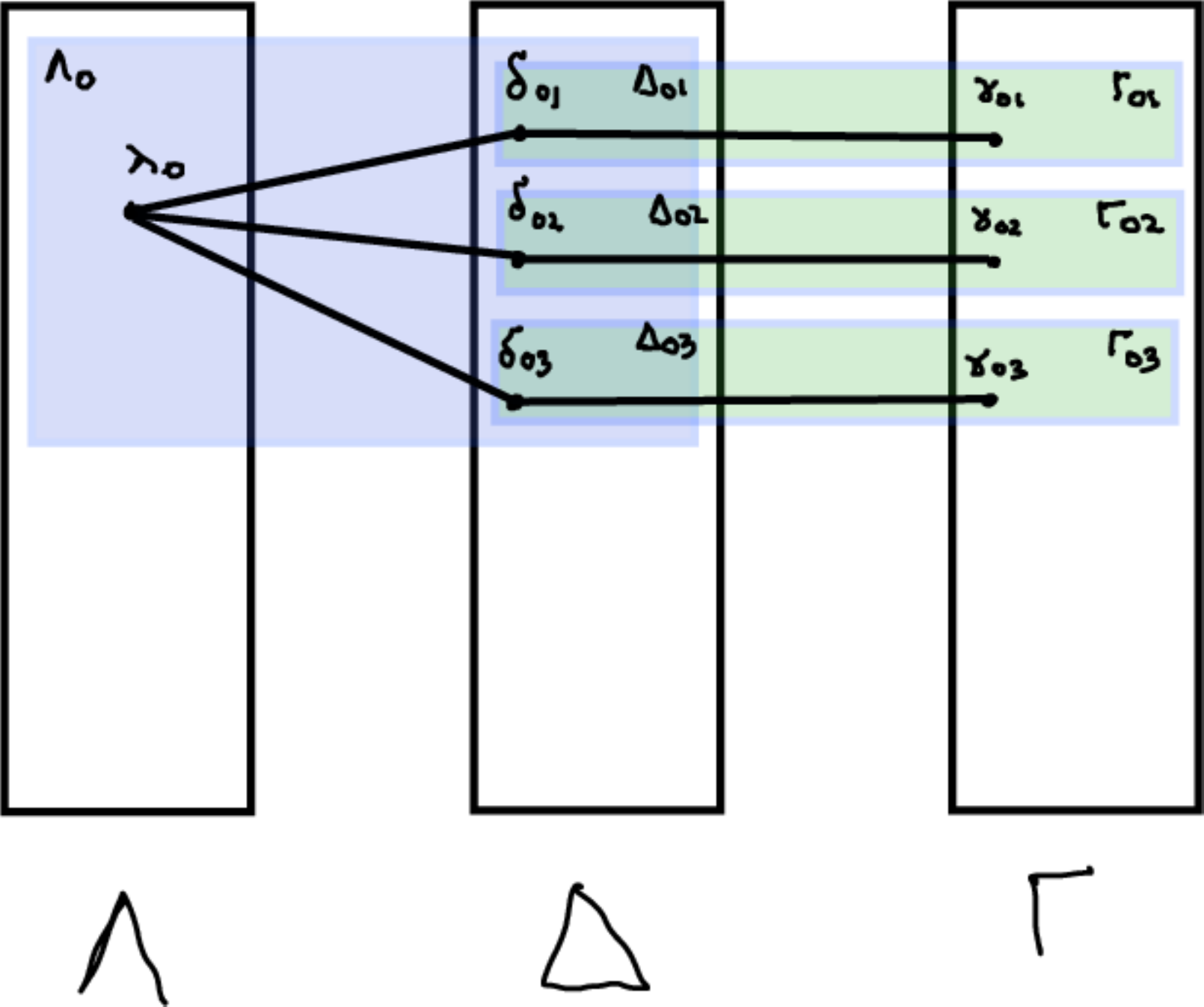}	
\caption{The vertices $\lambda_0,\delta_{01},\delta_{02},\delta_{03},\gamma_{01},\gamma_{02},\gamma_{03}$ in Lemma~\ref{270420a8}.\label{fig2}}
\end{figure}

Let $\Lambda_0\in\Lambda/e\sbs{\Lambda}{\Delta}$. By Theorem~\ref{080520p}(2) for $\Gamma=\Lambda$, there is a unique  $\Delta_0\in\Delta/e\sbs{\Delta}{\Lambda}$ such that $X_{\Delta_0\cup\Lambda_0}$ is complete bipartite. In view of inequality~\eqref{030520i}, $\Delta_0$ is the union of at least three pairwise distinct $\Delta_{01},\Delta_{02},\Delta_{03}\in\Delta/e\sbs{\Delta}{\Gamma}$. Again by Theorem~\ref{080520p}(2), there are uniquely determined $\Gamma_{01},\Gamma_{02},\Gamma_{03}\in\Gamma/e\sbs{\Gamma}{\Delta}$ such that $X_{\Delta_{0i}\cup\Gamma_{0i}}$ is complete bipartite, $i=1,2,3$. Choose arbitrarily the vertices
$$
\lambda_0\in \Lambda_0,\ \delta_{0i}\in\Delta_{0i},\ \gamma_{0i}\in\Gamma_{0i}\quad (i=1,2,3),
$$
the  subgraph $Y$ induced by these vertices is depicted in Fig.~\ref{fig2}.
It is easily seen that $Y$  is isomorphic to~$T_2$, a contradiction.\eprf\medskip

Let $\Delta,\Gamma\in F$. Without loss of generality we may assume that  $X_{\Delta\cup\Gamma}$ is not empty, for otherwise, $\nrm{e\sbs{\Delta}{\Gamma}}=1$. Next, denote by $\alpha$ the distinguished vertex of $X$. By the connectivity of $X$,  the number
$$
d=d(\alpha,\Delta\cup\Gamma)
$$ 
is a nonnegative integer.  By Lemma~\ref{regular}(2), this implies  that $\{d(\alpha,\Delta),d(\alpha,\Gamma)\}=\{d,d+1\}$. To prove inequality~\eqref{010520l}, we use induction on~$d$. By the symmetry between~$\Delta$ and~$\Gamma$, we may also assume that 
$$
d(\alpha,\Delta)=d\qaq d(\alpha,\Gamma)=d+1.
$$

When $d=0$, we have $\Delta=\{\alpha\}$. It follows that $X_{\Delta\cup\Gamma}$ (being nonempty) is complete bipartite with parts $\Delta$ and $\Gamma$. Therefore, $\nrm{e\sbs{\Delta}{\Gamma}}=1$, and we are done.\medskip

Let $d\ge 1$. Then there exists a vertex $\lambda$ of $X$ having a neighbor in $\Delta$ and such that $d(\alpha,\lambda)=d-1$. Denote by $\Lambda$ the fiber of $\cX$, containing~$\lambda$.  Then  $X_{\Lambda\cup\Delta}$ is nonempty and 
$d(\alpha,\Lambda\cup\Delta)=d-1$. By the induction hypothesis, this yields 
\qtnl{010520j}
\nrm{e\sbs{\Delta}\Lambda}\le 2. 
\eqtn

Assume on the contrary that  $\nrm{e\sbs{\Delta}{\Gamma}}\ge 3$. Then  $\nrm{e\sbs{\Delta}{\Gamma}}> \nrm{e\sbs{\Delta}\Lambda}$.  According to Theorem~\ref{080520p}(4), this shows that $e\sbs{\Delta}{\Gamma}\subsetneq e\sbs{\Delta}{\Lambda}$. Consequently, $\nrm{e\sbs{\Delta}{\Gamma}}=2\nrm{e\sbs{\Delta}{\Lambda}}$ by Lemma~\ref{270420a8}. In view of the assumption and  inequality~\eqref{010520j}, this is possible only if
\qtnl{020520a}
\nrm{e\sbs{\Delta}{\Gamma}}=4\qaq\nrm{e\sbs{\Delta}{\Lambda}}=2.
\eqtn

\begin{figure}[h]
	\includegraphics[scale=0.2]{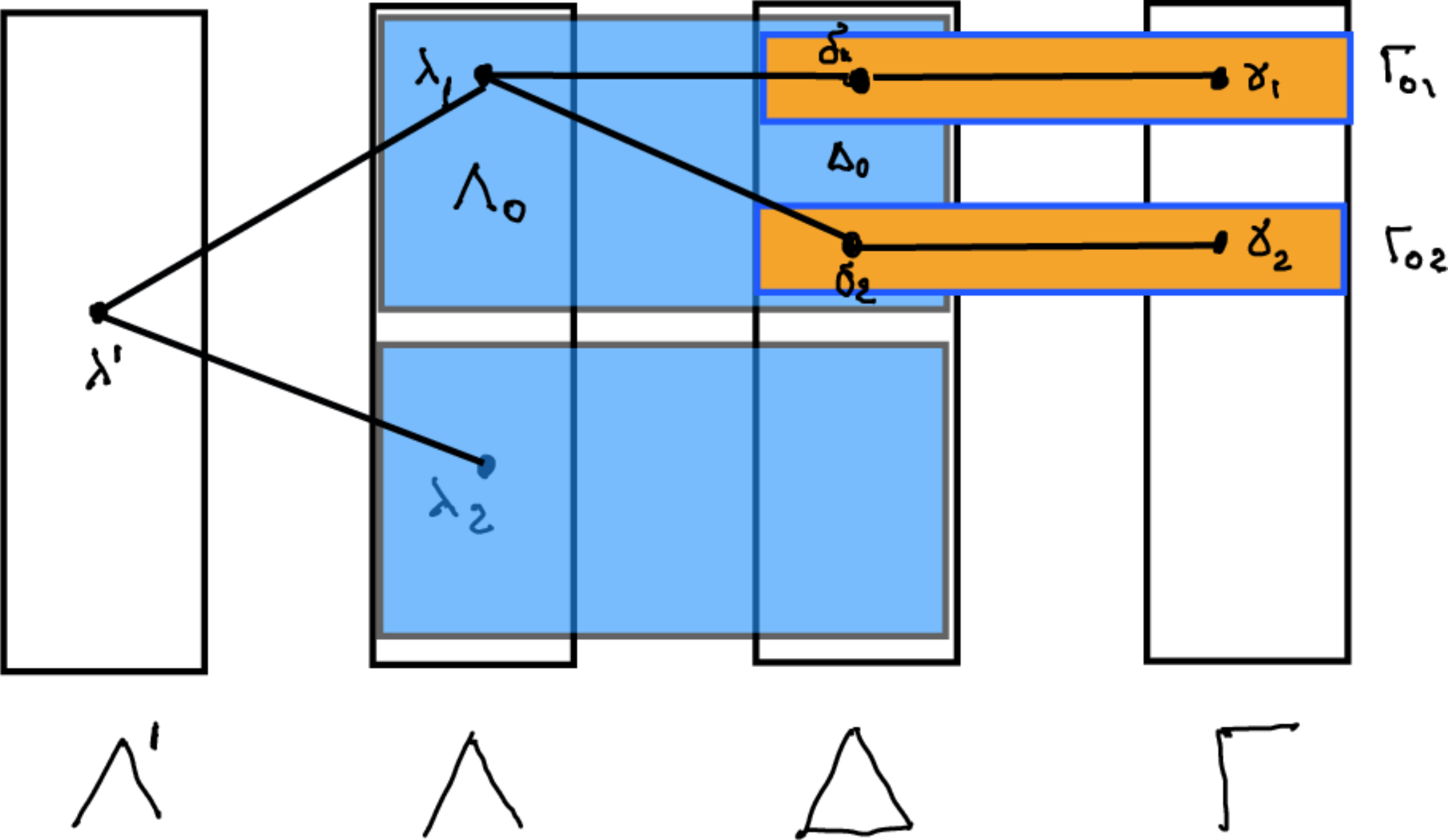}	
	\caption{The vertices $\lambda',\lambda_1,\lambda_2,\delta_1,\delta_2,\gamma_1,\gamma_2$ in Case~1.\label{fig3}}
\end{figure}

\begin{figure}[h]
\includegraphics[scale=0.2]{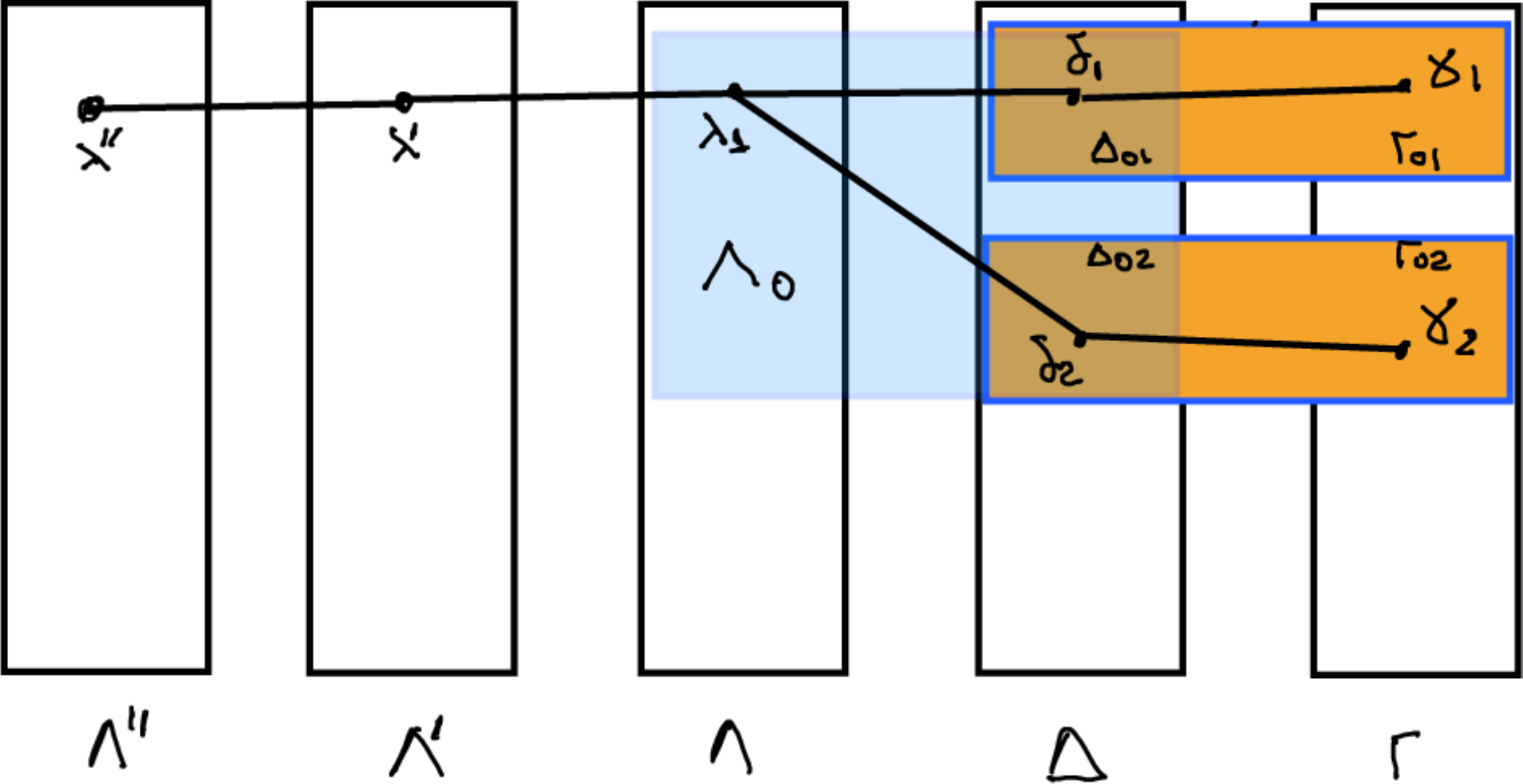}	
\caption{The vertices $\lambda'',\lambda',\lambda_1,\delta_1,\delta_2,\gamma_1,\gamma_2$ in Case~2.\label{fig4}}
\end{figure}

Let $\Lambda_0\in\Lambda/e\sbs{\Lambda}{\Delta}$. By Theorem~\ref{080520p}(2), for $\Gamma=\Lambda$, there is a unique $\Delta_0\in \Delta/e\sbs{\Delta}{\Lambda}$ such that $X_{\Delta_0\cup\Lambda_0}$ is complete bipartite. By formula~\eqref{020520a}, the class~$\Delta_0$ is the union of distinct $\Delta_{01},\Delta_{02}\in\Delta/e\sbs{\Delta}{\Gamma}$. Again by Theorem~\ref{080520p}(2), there are uniquely determined $\Gamma_{01},\Gamma_{02}\in\Gamma/e\sbs{\Gamma}{\Delta}$ such that  $X_{\Delta_{01}\cup\Gamma_{01}}$ and $X_{\Delta_{02}\cup\Gamma_{02}}$ are complete bipartite. Choose arbitrarily the vertices  
$$
\lambda_1\in\Lambda_0\quad 
\delta_1\in\Delta_{01},\ \delta_2\in\Delta_{02},\quad 
\gamma_1\in\Gamma_{01},\ \gamma_2\in\Gamma_{02}.
$$
Then the graph induced by these five  vertices in $X$ is isomorphic to a subgraph of the graphs depicted in Figs~\ref{fig3} and~\ref{fig4}. \medskip

By the second equality of formula~\eqref{020520a}, we have  $\Lambda\ne\{\alpha\}$ or, equivalently, $d(\alpha,\Lambda)\ge 1$. It follows that there is a vertex $\lambda'$ having a neighbor in $\Lambda$ and such that $d(\alpha,\lambda')=d(\alpha,\lambda)-1$. Denote by $\Lambda'$ the fiber of $\cX$, containing~$\lambda'$. We come to the final contradiction by considering two cases depending on whether or not the graph $X_{\Lambda'\cup\Lambda}$ is complete bipartite.\medskip

{\bf Case 1:} $X_{\Lambda'\cup\Lambda}$ is complete bipartite. Let us choose arbitrary vertex  $\lambda_2\in\Lambda\setminus\Lambda_0$; the obtained configuration is depicted in Fig~\ref{fig3}. 
By the assumption of the case, the  vertex $\lambda'$ is adjacent with $\lambda_1$ and $\lambda_2$. Furthermore, $\lambda'$ is adjacent with none of~$\gamma_1,\gamma_2$, because $d(\alpha,\delta')=d-2$, whereas $d(\alpha,\gamma_1)=(\alpha,\gamma_2)=d+1$. Thus the subgraph of~$X$, induced by the vertices $\lambda',\lambda_1,\lambda_2,,\delta_1,\delta_2,\gamma_1,\gamma_2$,  is isomorphic to~$T_2$, a contradiction.\medskip

{\bf Case 2:} $X_{\Lambda'\cup\Lambda}$  is not complete bipartite. Then $\Lambda'\ne\{\alpha\}$ or, equivalently, $d(\alpha,\Lambda')\ge 1$. It follows that there exists a vertex $\lambda''$ having a neighbor in $\Lambda'$ and such that $d(\alpha,\lambda'')=d(\alpha,\lambda')-1$. The obtained configuration is depicted in Fig~\ref{fig4}. 
As in the Case~1, the distance argument shows that the subgraph of $X$, induced by the vertices $\lambda'',\lambda'_1,\lambda_1,\delta_1,\delta_2,\gamma_1,\gamma_2$,   is isomorphic to~$T_2$, a contradiction.\eprf

\section{Proof of Theorems~\ref{main} and~\ref{main2}}\label{190520a}

{\bf Proof of Theorem~\ref{main2}.} Let $X$ be a connected   chordal bipartite $T_2$-free graph  with distinguished vertex, and $\cX=\WL(X)$. Denote by $\Omega$ the vertex set of~$X$ and by $t=t(\cX)$ the twin parabolic of~$\cX$. It suffices to verify that
\qtnl{270420a}
\nrm{t\sbo{\Delta}}\le 2\quad\text{for all}\ \,\Delta\in F,
\eqtn
where $F=F(\cX)$. Indeed, then the cardinality of  every fiber of the quotient coherent configuration $\cX_{\Omega/t}$ is at most~$2$. According to~\cite[Exercise~3.7.20]{CP}, this implies that $\cX_{\Omega/t}$ is separable. Thus, $\cX$ is separable by Lemma~\ref{160520a}.\medskip

To prove formula~\eqref{270420a}, let $\Delta\in F$. By Theorem~\ref{080520p}(4), there exists $\Gamma\in F$ such that the relation $e\sbs{\Delta}{\Gamma}$ is minimal possible. We claim that 
\qtnl{010520i}
e\sbs{\Delta}{\Gamma}\subseteq t\sbo{\Delta}.
\eqtn
In other words, we need to verify that  every vertices $\alpha$ and $\beta$ lying in the same class of the equivalence relation $e\sbs{\Delta}{\Gamma}$ are $\cX$-twins.  However, for each $\Lambda\in F$,
$$
(\alpha,\beta)\in e\sbs{\Delta}{\Gamma}\subseteq e\sbs{\Delta}{\Lambda}.
$$
It follows that the vertices $\alpha$ and $\beta$ are $X_{\Delta\cup\Lambda}$-twins for all~$\Lambda$. Consequently, they  are $X$-twins. Thus the claim follows from Lemma~\ref{twin1}.\medskip

To complete the proof, we note that $\nrm{e\sbs{\Delta}{\Gamma}}\le 2$ by Theorem~\ref{300420a}. Together with inclusion~\eqref{010520i}, this yields
$$
\nrm{t\sbo{\Delta}}\le \nrm{e\sbs{\Delta}{\Gamma}}\le 2,
$$
which completes the proof of formula~\eqref{270420a}.\eprf\medskip

{\bf Proof of Theorem~\ref{main}.} One can see that $K_{n,n}$ is a chordal bipartite $T_2$-free graph for all~$n$. Moreover, $\dimwl(K_{n,n})\ne 1$ for all $n\ge 3$, see \cite[Lemma~3.1(A)]{AKRV2017}. Thus  suffices to prove that $\dimwl(X)\le 3$ for every  chordal bipartite $T_2$-free graph~$X$.\medskip

By statement Lemma~\ref {290420b}(1), we may assume that $X$ is connected.  Take an arbitrary vertex  $\alpha$ of this graph. Then the coherent configuration $\WL(X_{\alpha})$ is separable by Theorem~\ref{main2}. By Lemma~\ref {290420b}(3), this shows that $\dimwl(X_\alpha)\le 2$. Thus using Lemma~\ref {290420b}(2), we obtain
$$
\dimwl(X)\le \dimwl(X_\alpha)+1\le 3,
$$
as required.\eprf

\end{document}